\documentclass[runningheads,a4paper]{llncs}

\usepackage{amsfonts,amssymb,amsmath}
\usepackage{natbib}
\usepackage{wrapfig}
\usepackage{llncsdoc}
\usepackage{epsfig,psfrag}
\usepackage{latexsym}
\usepackage{longtable}
\usepackage{tikz,pgf}
\usepackage{subfig}
\usepackage{wrapfig}

\usepackage{graphicx}
\usepackage{epstopdf}
\usepackage{color}

\newcommand{\bqn}{\begin{eqnarray}}
\newcommand{\eqn}{\end{eqnarray}}
\newcommand{\bq}{\begin{eqnarray*}}
\newcommand{\eq}{\end{eqnarray*}}

\usepackage{url}
\usepackage{hyperref}
\hypersetup{colorlinks,%
citecolor=black,%
filecolor=blue,%
linkcolor=red,%
urlcolor=blue,%
pdftex}

\newcommand{\blue}[1]{{\color{blue} #1}}

\begin{document}

\title{Diffusion Equations for Medical Images}
 \titlerunning{Diffusion Equations for Medical Images}
 
\author{Moo K. Chung
 }
\institute{
University of Wisconsin-Madison, USA\\
\vspace{0.3cm}
\blue{\tt mkchung@wisc.edu}
}
\authorrunning{Chung}

\maketitle



\pagenumbering{arabic}

\index{Diffusions}

In brain imaging, the image acquisition and processing processes themselves are likely to introduce noise to the images. It is therefore imperative to reduce the noise while preserving the geometric details of the anatomical structures for various applications. Traditionally Gaussian kernel smoothing has been often used in brain image processing and analysis. However, the direct application of Gaussian kernel smoothing tend to cause various numerical issues in irregular domains with boundaries. For example, if one uses large bandwidth in kernel smoothing in a cortical bounded region, the smoothing will blur signals across boundaries. So in kernel smoothing and regression literature, various ad-hoc procedures were introduce to remedy the boundary effect.

Motivated by \citet{perona.1990}, diffusion equations have been widely used in brain imaging as a form of noise reduction.  
The most natural straightforward way to smooth images in irregular domains with boundaries is to formulate the problem as boundary value problems using  partial differential equations. Numerous diffusion-based techniques have been developed in  image processing \citep{sochen.1998,malladi.2002,tang.1999,taubin.2000,andrade.2001,chung.2001.diffusion,chung.2003.ni,chung.2005.IPMI,chung.2004.ISBI,cachia.MIA.2003,cachia.TMI.2003,joshi.2009}. In this paper, we will overview the basics of isotropic diffusion equations and explain how to solve them on regular grids and irregular grids such as graphs.

\section{Diffusion as a Cauchy problem}
\index{diffusion}
\index{diffusion!Cauchy problem}

Consider $\mathcal{M} \in \mathbb{R}^d$ to be a compact
differentiable manifold. Let $L^2(\mathcal{M})$ be the space of square integrable functions
in $\mathcal{M}$ with inner product \bqn \langle g_1, g_2 \rangle =
\int_{\mathcal{M}} g_1(p)g_2(p) \;d\mu(p),\label{eq:inner}\eqn where
$\mu$ is the Lebegue measure such that $\mu(\mathcal{M})$ is the
total volume of $\mathcal{M}$. The norm $\| \cdot \|$ is
defined as $$\|g\| = \langle g,g \rangle^{1/2}.$$
The linear partial differential operator $\mathcal{L}$ is {\em
self-adjoint} if
$$\langle g_1, \mathcal{L}g_2 \rangle = \langle
\mathcal{L}g_1, g_2 \rangle$$ for all $g_1, g_2 \in
L^2(\mathcal{M})$. Then the eigenvalues $\lambda_j$ and
eigenfunctions $\psi_j$ of the operator $\mathcal{L}$ are obtained
by solving \bqn \mathcal{L} \psi_j = \lambda_j
\psi_j.\label{eq:eigen}\eqn 
Often (\ref{eq:eigen}) is written as
\bq \mathcal{L} \psi_j = -\lambda_j
\psi_j\eq
so care should be taken in assigning the sign of eigenvalues.

\begin{theorem}
\label{lem:eigen}
The eigenfunctions $\psi_j$ are orthonormal.
\end{theorem}
{\em Proof.}
Note $\langle \psi_i, \mathcal{L} \psi_j \rangle = \lambda_j \langle \psi_i, \psi_j \rangle$. On the other hand, 
$\langle \mathcal{L} \psi_i, \psi_j \rangle = \lambda_i \langle \psi_i, \psi_j \rangle$. Thus
$$ (\lambda_i - \lambda_j )\langle \psi_i, \psi_j \rangle = 0.$$
For any $\lambda_i \neq \lambda_j$, $\langle \psi_i, \psi_j \rangle=0$, orthogonal. 
For $\psi_j$ to be orthonormal, we need $\langle \psi_j, \psi_j \rangle =1$. This is simply done by 
absorbing the constant multiple into $\psi_j$. Thus, $\{ \psi_j \}$ is orthonormal. $\square$

In fact $\psi_j$ is the basis in $L^2(\mathcal{M})$. Consider 1D eigenfunction problem 
$$\frac{\partial^2}{\partial x^2} \psi_j (x) = - \lambda_j \psi_j(x)$$
in interval $[-l, l]$. We can easily check that 
$$\psi_{1j} = \cos \big( \frac{j \pi x}{l} \big), \; \psi_{2j} = \sin \big( \frac{j \pi x}{l} \big), \; j = 1,2, \cdots$$
are eigenfunctions corresponding to eigenvalue $\lambda_j = \big( \frac{j \pi }{l} \big)^2$.
Also $\psi_{10} = 1$ is trivial first eigenfunction corresponding to $\lambda_0 =0$.  
The multiplicity of eigenfunctions is caused by the symmetric of interval $[-l, l]$. Based on trigonometric  formula, we can show that the eigenfunctions are orthogonal.
\bq \int_{-l}^l \psi_{1i}(x) \psi_{1j}(x) \; dx &=& 0 \mbox{ if } i \neq j \\
\int_{-l}^l \psi_{2i}(x) \psi_{2j}(x) \; dx &=& 0 \mbox{ if } i \neq j\\
\int_{-l}^l \psi_{1i}(x) \psi_{2j}(x) \; dx &=& 0 \mbox{ for any } i,j 
\eq
From $\psi_{1j}^2 (x) + \psi_{2j}^2 (x) =1$ and due to symmetry
$$\int_{-l}^l \psi_{1j}^2 (x)  \; dx =  \int_{-l}^l \psi_{2j}^2 (x) \; dx = l.$$
Thus 
\bq\psi_{10} &=& \frac{1}{\sqrt{2l}}, \\
\psi_{1j} &=& \frac{1}{\sqrt{l}} \cos \big( \frac{j \pi x}{l} \big), \\
\psi_{2j} &=& \frac{1}{\sqrt{l}}\sin \big( \frac{j \pi x}{l} \big), \; j =1, 2, \cdots
\eq
are orthonormal basis in $[-l,l]$.

Consider a Cauchy problem of the following form: 
\bqn 
\frac{\partial g}{\partial t}(p,t) +
\mathcal{L}g(p,t) =0 \label{eq:cauchy1}, g(p,t=0) =f(p),\eqn
where $t$ is time variable and $p$ is spatial variable. 

The initial functional data $f(p)$ can be further
stochastically modeled as
\bqn f(p) = \nu(p) + \epsilon(p),\label{eq:stochastic1}\eqn where
$\epsilon$ is a stochastic noise modeled as a zero-mean Gaussian
random field, i.e., $\mathbb{E} \epsilon(p)=0$ at each point $p$ and $\nu$ is the unknown signal to be estiamted. PDE
(\ref{eq:cauchy1}) diffuses noisy initial data $f$ over time and estimate
the unknown signal $\nu$ as a solution. Diffusion time $t$ controls the
amount of smoothing and will be termed as the {\em bandwidth}. The unique solution to
equation (\ref{eq:cauchy1}) is given as follows. This is a heuristic proof and more rigorous proof is given later. \\

\begin{theorem}\label{thm:cauchy} For the self-adjoint linear differential operator
$\mathcal{L}$, the unique solution of the Cauchy problem
\bqn 
\frac{\partial g}{\partial t}(p,t) +
\mathcal{L}g(p,t) =0 , g(p,t=0) =f(p) \label{eq:cauchy2}\eqn
is given by
\bqn g(p,t) = \sum_{j=0}^{\infty} e^{-\lambda_j t} \langle f, \psi_j
\rangle \psi_j(p).\label{eq:cauchy.solution}\eqn
\end{theorem}
{\em Proof.} For each fixed $t$, since $g \in L^2(\mathcal{M})$,
$g$ has expansion \bqn g(p,t) = \sum_{j=0}^{\infty} c_j(t) \psi_j(p)
\label{eq:expansion}.\eqn Substitute equation (\ref{eq:expansion})
into (\ref{eq:cauchy2}). Then we obtain \bqn \frac{\partial}{\partial t} c_j(t) +
\lambda_j c_j(t)=0\label{eq:ODE}\eqn 
for all $j$. The solution of equation
(\ref{eq:ODE}) is given by $$c_j(t) = b_je^{-\lambda_j t}.$$ So we
have solution $$g(p,t) = \sum_{j=0}^{\infty} b_j e^{-\lambda_j t}
\psi_j(p).$$ At $t=0$, we have
$$g(p,0) = \sum_{j=0}^{\infty} b_j \psi_j(p) = f(p).$$
The coefficients $b_j$ must be the Fourier coefficients $\langle
f,\psi_j \rangle$ and they are uniquely determinded. $\square$

\begin{figure}[t]
\centering
\includegraphics[width=1\linewidth]{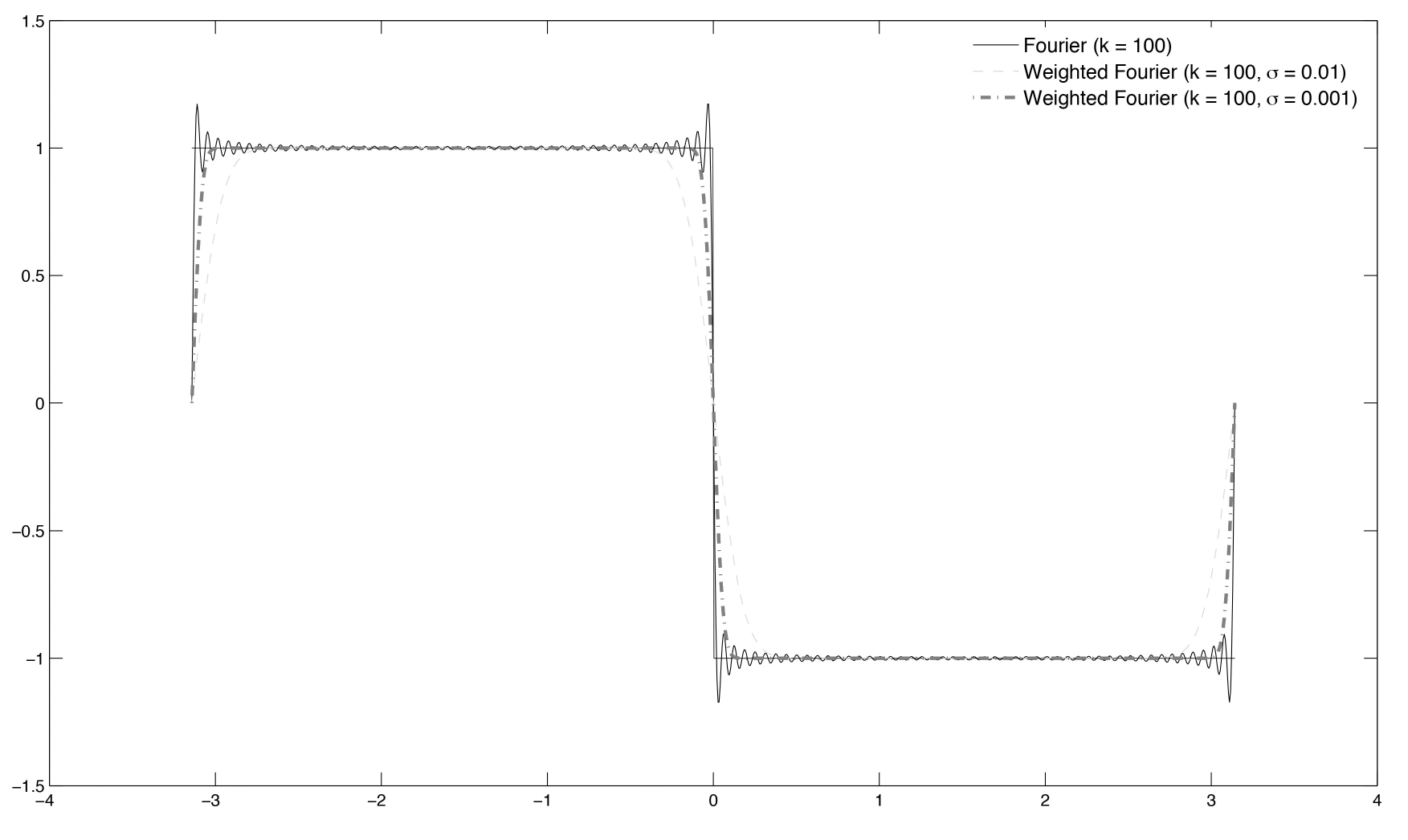}
 \caption{Reduction of Gibbs phenomenon in weighted Fourier series (WFS). The figure was generated by Yuan Wang of University of South Carolinia \citep{wang.2018.annals}.}
\label{fig:Lhkernel}
 \end{figure}
 
The implication of Theorem \ref{thm:cauchy} is obvious. The
solution decreases exponentially as time $t$ increases and smoothes
out high spatial frequency noise much faster than low
frequency noise. This is the basis of many
of PDE-based image smoothing methods. PDE involving self-adjoint linear
partial differential operators such as the Laplace-Beltrami operator
or iterated Laplacian have been widely used in medical image
analysis as a way to smooth either scalar or vector data along
anatomical boundaries \citep{andrade.2001,bulow.2004,chung.2003.ni}. These methods directly solve PDE using standard numerical techniques such as the finite
difference method (FDM) or the finite element method (FEM). The main shortcoming of 
solving PDE using FDM or FEM is the numerical instability and the
complexity of setting up the numerical scheme. The analytic approach called weighted Fourier series (WFS) differs from these previous methods in such a way that we only need
to estimate the Fourier coefficients in a hierarchical fashion to solve PDE \citep{chung.2007.tmi}.

\begin{example}Consider 1D differential operator $\mathcal{L} = \frac{\partial^2}{\partial x^2}$. The corresponding Cauchy problem is 1D diffusion equation 
$$ \frac{\partial g} {\partial t} (p,t) + \frac{\partial^2 g}{\partial x^2} (p,t) = 0, g(p,t=0)=f(p), p \in [-l,l].$$
Then the solution of this problem is given by Theorem \ref{thm:cauchy}:
$$g(p,t) = a_0 \psi_{10} + \sum_{j=1}^{\infty} a_j e^{-\lambda_j t} \psi_{1j}(p)+ b_j e^{-\lambda_j t}\psi_{2j}(p),$$
where
\bq a_0 &=& \frac{1}{\sqrt{2l}} \int_{-l}^l f(p) \; dp, \\
a_j &=& \frac{1}{\sqrt{l}} \int_{-l}^l f(p)  \cos \big( \frac{j \pi x}{l} \big) \; dp, \\
b_j &=& \frac{1}{\sqrt{l}}   \int_{-l}^l f(p)\sin  \big(  \frac{j \pi x}{l} \big) \; dp
\eq
for $j =1, 2, \cdots.$
\end{example}

\section{Finite difference method}
\index{finite difference}

One way of solving diffusion equations numerically in  to use finite differences. We will discuss how to differentiate images. There are numerous techniques for differentiation proposed in literature. We start with simple example of image differentiation in 2D image slices. Consider image intensity $f(x,y)$ defined on a regular grid, i.e., $(x,y) \in \mathbb{Z}^2$. Assume the pixel size is $\delta x$ and $\delta y$ in $x$- and $y$-directions. The partial derivative along the $x$-direction of image $f$ is approximated by the finite difference:

\bq \frac{\partial f}{\partial
x}(x,y) = \frac{f(x + \delta x, y) - f(x, y)}{\delta x}.
\label{eq:diff1}
\eq
The partial derivative along the $y$-direction of image $f$ is approximated similarly. $\frac{\partial f}{\partial
x}(x,y)$ and $\frac{\partial f}{\partial
y}(x,y)$  are called the {\em first order derivatives}. Then the {\em second order derivatives} are defined by taking the finite difference twice:

\bq
\frac{\partial^2 f}{\partial x^2}(x,y) &=&
\Big[ \frac{f(x + \delta x, y) - f(x,y)}{\delta x} -
\frac{f(x, y) - f(x-\delta x,y)}{\delta x}   \Big] /\delta x  \\
&=& \frac{f(x+ \delta x, y) - 2f(x,y) + f(x-\delta x, y)}{\delta x^2}
\label{eq:Lx}
 \eq

Similarly, we also have
 \bq
 \frac{\partial^2 f}{\partial y^2}(x,y) &=& \frac{f(x, y + \delta y) - 2f(x,y) + f(x, y - \delta y)}{\delta x^2}.
\label{eq:Ly}
 \eq
Other partial derivatives such as $\frac{\partial^2 f}{\partial x \partial y}$ are  computed similarly.

\subsection{1D diffusion by finite difference}

\begin{figure}[t!]
\centering
\includegraphics[width=1\linewidth]{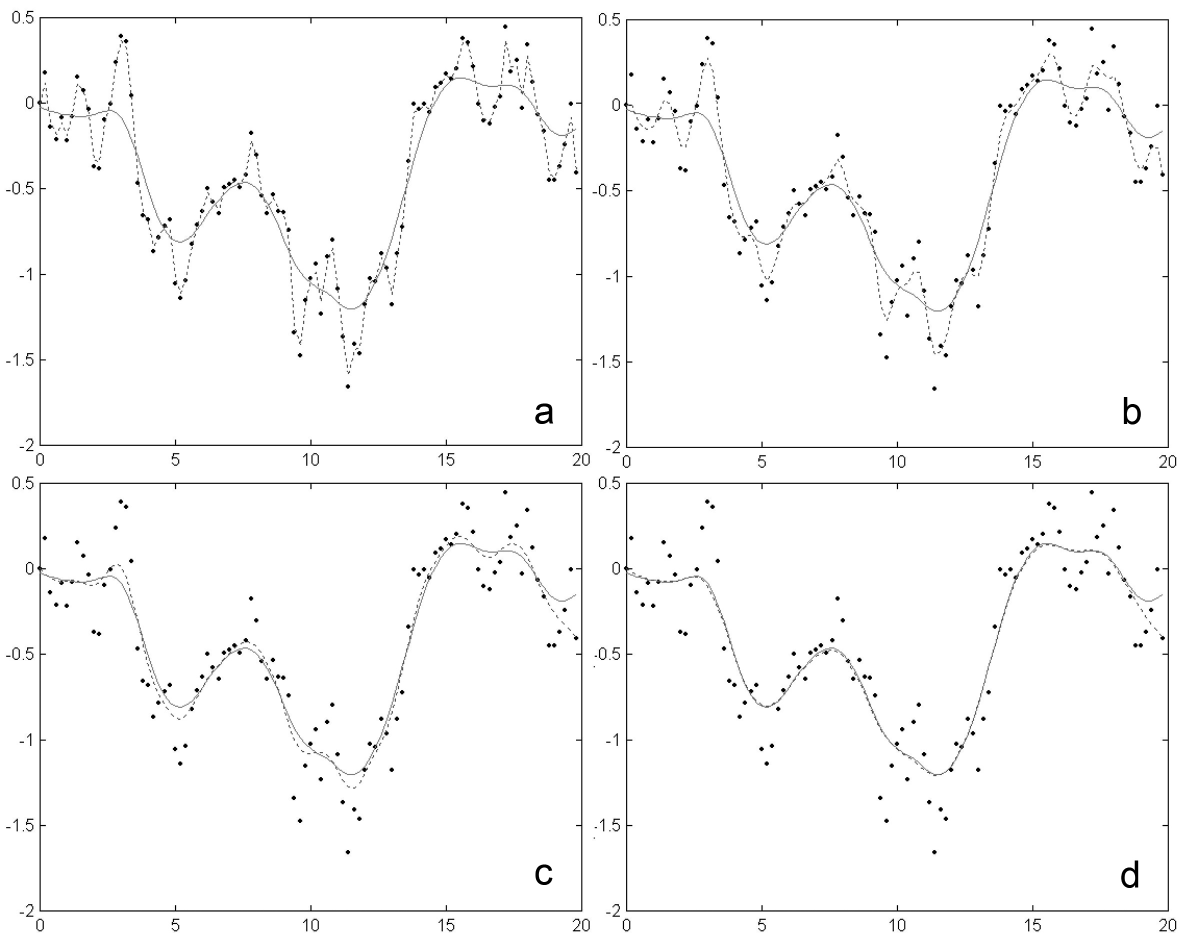}
 \caption{Gaussian kernel smoothing  (solid line) and diffusion smoothing (dotted  line). (a) before diffusion (b) after 0.05 seconds (5 iterations) (c) after 0.25 seconds (25 iterations) (d) after 0.5 seconds (50 iterations).}
\label{fig:diffusion-1Dsmoothing}
 \end{figure}
 
Let us implement 1D version of diffusion equations (Figure \ref{fig:diffusion-1Dsmoothing}). Suppose we have a smooth function $f(x,t)$ which is a function of position $x \in \mathbb{R}$ and time $t \in \mathbb{R}^+$. 1D isotropic heat equation is then defined as
\bqn \frac{\partial f}{\partial t} = \frac{d^2 f}{d
x^2}
\label{eq:1Ddiffs}
\eqn
with initial condition $f(x, t=0) = g(x)$. Differential equation (\ref{eq:1Ddiffs}) is then discretized as
\bqn f(x, t + \delta t) &=&  f(x,t)  + {\delta t} \frac{d^2 f}{d x^2}(x,y).
\label{eq:2Diffs} \eqn
With $t_k = k \delta t$ and starting from $t=0$, (\ref{eq:2Diffs}) can be written as
\bqn f(x, t_{k+1}) = f(x, t_k) + 
 \delta t \frac{f(x+ \delta x,t_k) - 2f(x,t_k) + f(x-\delta x,t_k)}{\delta x^2}. \label{eq:1diffusionFD}
\eqn

The above finite difference gives the solution at time $t_{k+1}$. To obtain the solution at any time, it is necessary to keep iterating many times with
very small $\delta t$. If $\delta t$ is too small, the computation
is slow. If it is too large, the finite difference will diverge. Then the problem is finding the largest $\delta t$ that
grantee the convergence.


 \begin{figure}[t]
\centering
\includegraphics[width=0.8\linewidth]{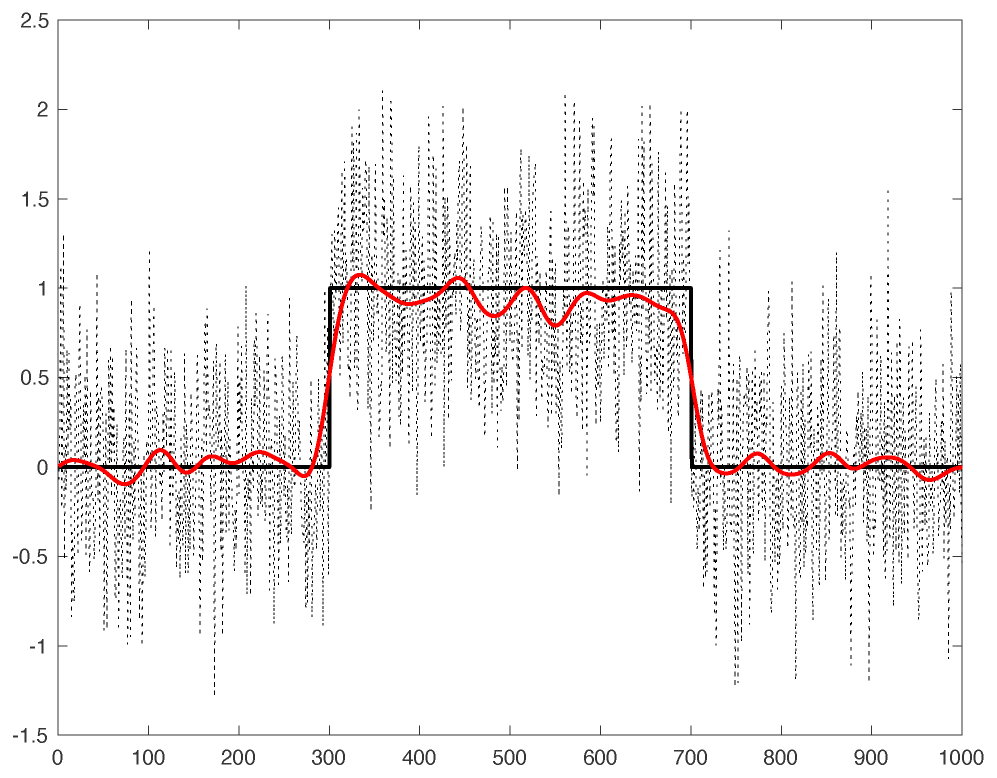}
 \caption{Diffusion of simulated data (dotted line) and the ground truth (black line). The both methods using {\tt conv.m} and {\tt toeplitz.m} all converge to the red line.}
\label{fig:1Ddiffusion}
 \end{figure}

Numerically (\ref{eq:1diffusionFD}) is solved in MATLAB follows. We start with generating a step function as the ground truth (black line in Figure \ref{fig:1Ddiffusion}). We then add $N(0,0.5^2)$ noise.

\begin{verbatim}
x=1:1000;
noise=normrnd(0, 0.5, 1,1000);
signal= [zeros(1,300) ones(1,400) zeros(1,300)];
figure; plot(signal, 'k', 'LineWidth',2);
y= signal + noise; 
hold on; plot(y, ':k');
\end{verbatim}

The 2nd order finite difference is coded as {\tt L=[1 -2 1]}. Then it is convoluted with 3 consecutive data at a time in the code blow.

\begin{verbatim}
L = [1 -2 1]
g=y; 
for i=1:10000  
    Lg = conv(g,L,'same');
    g = g+ 0.01*Lg;
end
hold on;plot(g, 'b', 'LineWidth', 2);   
\end{verbatim}

Since the Laplacian is a linear operator, the above convolution can be written as the matrix multiplication. Any linear operation can be discretely encoded as matrix multiplication. Here the Laplacian is encoded using a Toeplitz matrix:

\begin{verbatim}
c=zeros(1,1000);
c(1:2)=[-2 1];
r=zeros(1,1000);
r(1:2)=[-2 1];
L = toeplitz(c,r);
\end{verbatim}

The first 5 columns and rows of the Toelitz matrix {\tt L} is given by
\begin{verbatim}
 L(1:5,1:5)
 ans =
     -2     1     0     0     0
      1    -2     1     0     0
      0     1    -2     1     0
      0     0     1    -2     1
      0     0     0     1    -2
\end{verbatim}

The diffusion is solved by sequential summation of matrix multiplications
\begin{verbatim}
g=y'; 
for i=1:10000  
    g = g+ 0.01*L*g;
end
hold on;plot(g, 'g', 'LineWidth', 2);     
\end{verbatim}
The first and last rows of the the Toelitz matrix {\tt L} is {\tt [-2 1 0]} and {\tt [0 1 -2]}, which is different from the 2nd order finite difference {\tt [1 -2 1]}. It will not matter since that may be viewed as the discrete Laplacian in the boundary.  The matrix form of Laplacian can be used in conjunction with the recently developed polynomial approximation method for solving heat diffusion on manifolds \citep{huang.2020.TMI}. \\

{\em Discrete maximum principle.} Since the diffusion
smoothing and kernel smoothing are equivalent, The diffused signal
$f(x,t_{k+1})$ must be bounded by the minimum and the maximum of
signal \citep{chung.2003.cvpr}. Let $x_{i-1}, x_i, x_{i+1}$ be some points with gap $\delta x$.

\bq f(x_i,t_{k+1}) &=& f(x_i,t_k) + \delta t
\frac{d^2 f}{d x^2}(x_i,t_k)\\
&\leq& \max \big[f(x_{i-1},t_k), f(x_i,t_k),f(x_{i+1},t_k)
\big].\eq
Similarly, we can bound it below. Thus, the time step should be bounded by

$$\delta t \leq \max \Big[ \Big|\frac{f(x_{i-1},t_j) - f(x_i,t_j)}
{\frac{d^2 f}{d
x^2}}\Big|,\Big|\frac{f(x_{i+1},t_j)-f(x_i,t_j)}{\frac{d^2
f}{d x^2}}\Big|\Big].$$

\subsection{Diffusion in $n$-dimensional grid}
\index{diffusion}

In 2D,  let $(x_i,y_i)$ be pixels around $(x,y)$ including $(x,y)$ itself.
Then using the 4-neighbor scheme, Laplcian of $f(x,y)$ can be written as
$$\Delta f(x,y) = \sum_{i,j} w_{ij}f(x_i,y_i),$$
where the Laplacian matrix is  given by
$$(w_{ij}) = \left(%
\begin{array}{ccc}
  0 & 1 & 0 \\
  1 & -4 & 1 \\
  0 & 1 & 0 \\
\end{array}%
\right).
$$
Note that $\sum_{ij} w_{ij}=1$. 

Extending it further, consider $n$D. Let $x=(x_1, x_2, \cdots, x_n)$ be the coordinates in $\mathbb{R}^n$. 
Laplacian $\Delta$ in $\mathbb{R}^n$ is defined as
$$\Delta f = \frac{\partial^2 f}{\partial x_1^2} + \cdots + \frac{\partial^2 f}{\partial x_n^2}.$$
Assume we have a $n$-dimensional hyper-cube grid of size is 1. Then we have
\bq   \Delta f(x) &=& f(x_1 \pm 1,\cdots, x_n)  + \cdots + f(x_1, \cdots, x_n \pm 1)  \\
&&- 2n f(x,y).
\eq
This uses $2n$ closest neighbors of voxel $x$  to approximate the Laplacian.

It is also possible to incorporate $2^n$ corners $(x_1 \pm 1, \cdots, x_n \pm 1)$ along with the $2n$ closest neighbors for a better approximation of the Laplacian. In particular in 2D,  we can obtain a more accurate finite difference formula for 8-neighbor Laplacian:

$$(w_{ij}) = \frac{1}{9}\left(%
\begin{array}{ccc}
  1 & 1 & 1 \\
  1 & -8 & 1 \\
  1 & 1 & 1 \\
\end{array}%
\right).$$

Based on the estimation Laplacian on discrete grid, diffusion equation 
\bq \frac{\partial f}{\partial t} = \Delta f
\eq
is discretized as
\bqn f(x, t_{k+1}) = f(x, t_k) + \delta t \sum_{i,j} w_{ij}f(x_i, y_i). \label{eq:1Dfinite}
\eqn
with $t_k = k \delta t$ and starting from $t_1=0$. From (\ref{eq:1Dfinite}), we can see that the diffusion equation is solved by iteratively applying 
convolution with weights $w_{ij}$. In fact, it can be shown that the solution of diffusion is given by kernel smoothing.

%


%

\section{Laplacian on planner graphs}
\index{Laplacian}
\index{Laplacian!planar graphs}

In a previous section, we showed
how to estimate the Laplacian in a regular grid. Now we
show how to estimate Laplacian in irregular grid such as graphs and polygonal
surfaces in $\mathbb{R}^2$. The question is how one estimate Laplacian
or any other differential operators on a graph.
Assume we have observations $Y_i$ at each point $p_i$, which is
assumed to follow additive model 
$$Y_i = \mu(p_i) + \epsilon(p_i), \;
p_i \in \mathbb{R}^2$$ 
where $\mu$ is a smooth continuous function
and $\epsilon$ is a zero mean Gaussian random field. We want to
estimate at some node $p_i$ on a graph: 
$$\Delta \mu (p_0)= \frac{\partial^2 \mu}{\partial x^2}\Big|_{p_0} +
\frac{\partial^2 \mu}{\partial y^2}\Big|_{p_0}.$$ Unfortunately,
the geometry of the graph forbid direct application of
finite difference scheme. To answer this problem, one requires the finite element method (FEM) \citep{chung.2001.thesis}. However, we can use  a more elementary technique called {\em polynomial regression}.

Let $p_i = (x_i,y_i)$ be 
the coordinates of the vertices of the graph or polygonal surface. 
Let $p_i$ be the neighboring vertices of $p_0$. We estimate the Laplacian
at $p_0$ by fitting a quadratic polynomial of the form
 \bqn \mu(u,v) =\beta_0 + \beta_1u + \beta_2v + \beta_3u^2 +
\beta_4uv + \beta_5v^2.\label{eq:quadratic} \eqn
We are basically assuming the unknown signal $\mu$ to be the quadratic form (\ref{eq:quadratic}). 
Then the parameters $\beta_i$ are estimated by solving the normal
equation: 
\bqn Y_i =\beta_0 + \beta_1x_i + \beta_2y_i +
\beta_3x_i^2 + \beta_4x_iy_i + \beta_5y_i^2\eqn
for all $p_i$ that is neighboring $p_0$. For simplicity, we may assume $p_0$ is translated to the origin, i.e.,
$x_0=0, y_0 =0$.

Let $Y=(Y_1,\cdots,Y_m)^\top$, $\beta = (\beta_0,\cdots,\beta_5)^\top$ and 
design matrix
$$\mathbb{X}= \begin{pmatrix}
1& x_1 & y_1 & x_1^2 & x_1y_1 & y_1^2\\
1& x_2 & y_2 & x_2^2 & x_2y_2 & y_2^2\\
  \hdotsfor[2]{6} \\
1&   x_m & y_m & x_m^2 & x_my_m & y_m^2
\end{pmatrix}.$$
Then we have the following matrix equation
$$Y =\mathbb{X}\beta.$$ The unknown coefficients vector $\beta$ is
estimated by the usual least-squares method:
$$\widehat \beta = (\widehat\beta_0,\dots,\widehat \beta_5)^\top
=(\mathbb{X}^\top \mathbb{X})^{-}\mathbb{X}^\top Y,$$ where $^-$ denotes
generalized inverse, which can be obtained through the singular
value decomposition (SVD). Note that $\mathbb{X}^\top\mathbb{X}$ is
nonsingular if $m <6$. In Matlab, {\tt pinv} can be used to compute the generalized inverse, which is often called the {\em pseudo inverse}.\\

The generalized inverse often used is that of Moore-Penrose. 
It is usually defined as matrix $\mathbb{X}$ satisfying four conditions
$$\mathbb{X}\mathbb{X}^{-}\mathbb{X}=\mathbb{X}, \; \mathbb{X}^{-}\mathbb{X}\mathbb{X}^{-}
 = \mathbb{X}^{-},$$
$$(\mathbb{X}\mathbb{X}^{-})^\top=\mathbb{X}\mathbb{X}^{-}, \; (\mathbb{X}^{-}\mathbb{X})^\top=
\mathbb{X}^{-}\mathbb{X}.$$
Let $\mathbb{X}$ be $m \times p$ matrix with $m \geq p$. Then SVD
of $\mathbb{X}$ is
$$\mathbb{X}=UDV^\top,$$
where $U_{m \times p}$ has orthonormal columns, $V_{p \times p}$ is orthogonal, and
$D_{p \times p}=Diag(d_1,\cdots,d_p)$ is diagonal with non-negative elements
and . Let
$$D^{-}=Diag(d_1^{-},\cdots,d_p^{-}),$$
where $d_i^{-} = 1/d_i$ if $d_i
\neq 0$ and $d_i^{-} = 0$ if $d_i=0$. Then it can be shown that
the Moore-Penrose generalized inverse is given by
$$\mathbb{X}^{-}=VD^{-}U^\top.$$

Once we estimated the parameter vector $\beta$,  the Laplacian of is
$$\Delta \mu (p_0)= 2\widehat \beta_3 + 2\widehat
\beta_5.$$ 


\section{Graph Laplacian}
\index{graph!Laplacian}

Now we generalize volumetric Laplacian in previous sections to graphs. Let $G = (V,E)$ be a graph with node set $V$ and edge set $E$. We will simply index the node set as $V = \{1,2 ,\cdots, p\}$. If two nodes $i$ and $j$ form an edge, we denote it as $i \sim j$. Let $W=(w_{ij})$ be the edge wight. The adjacency matrix of $G$ is often used as the edge weight. Various forms of graph Laplacian have been proposed \citep{chung.1997} but the most often used standard form $L=(l_{ij})$ is given by 
$$l_{ij} =\left(
\begin{array}{cc}
  -w_{ij}, &  i \sim j\\
  \sum_{i \neq j} w_{ij}, &  i = j \\
  0, & \mbox{ otherwise}      
\end{array}\right .$$
Often it is defined with the sign reversed such that
$$l_{ij} =\left(
\begin{array}{cc}
  w_{ij}, &  i \sim j\\
  -\sum_{i \neq j} w_{ij}, &  i = j \\
  0, & \mbox{ otherwise}      
\end{array}\right .$$

The graph Laplacian $L$ can then be written as $$L = D- W,$$ where $D=(d_{ij})$ is the diagonal matrix with $d_{ii} = \sum_{j =1}^n w_{ij}$.
Here, we will simply use the adjacency matrix so that the edge weights $w_{ij}$ are either 0 or 1. In Matlab, Laplacian {\tt L} is simply computed from the adjacency matrix {\tt adj}:

\begin{verbatim}
n=size(adj,1); 
adjsparse = sparse(n,n);
adjsparse(find(adj))=1;
L=sparse(n,n);
GL = inline('diag(sum(W))-W'); 
L = GL(adjsparse);
\end{verbatim}

We use the sparse matrix format to reduce the memory burden for large-scale computation.

\begin{figure}[t]
\centering
\includegraphics[width=1\linewidth]{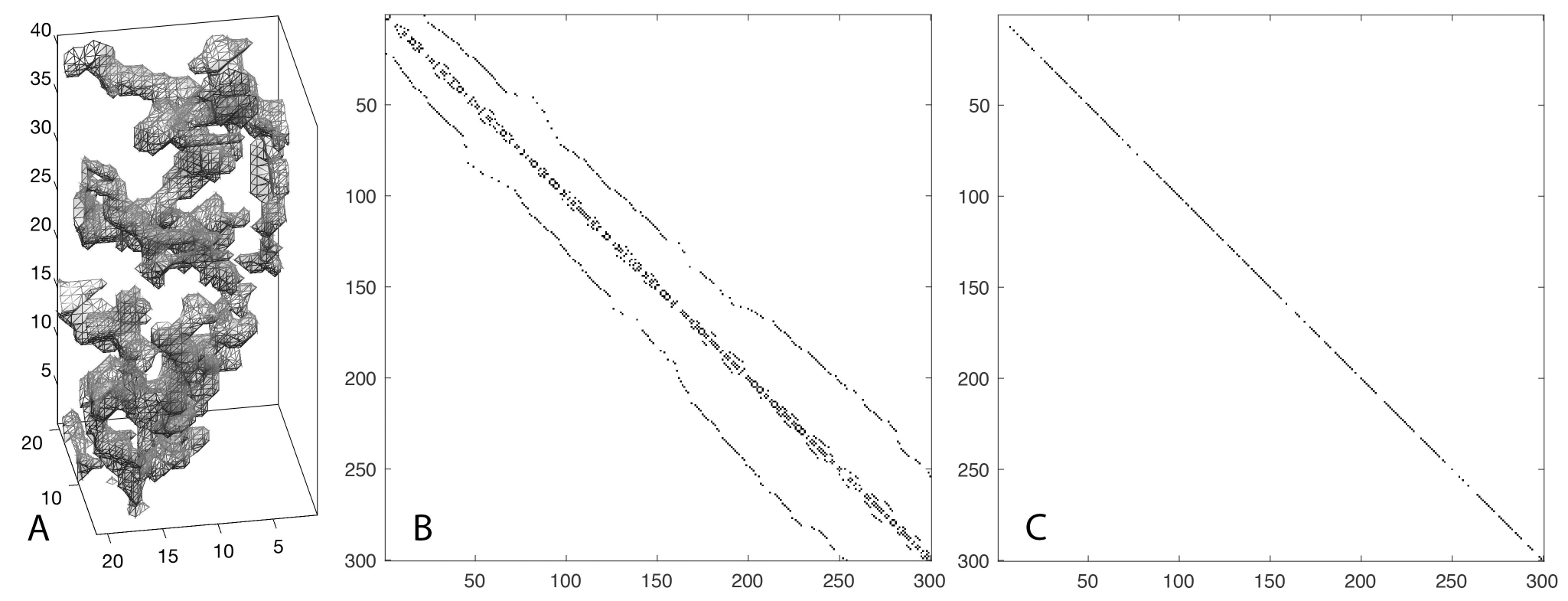}
\caption{A. Part of lung blood vessel obtained from CT. B. Adjacency matrix obtained from 4-neighbor connectivity. C. Laplace matrix obtained from the adjacency matrix.}
\label{fig:Laplacian-adj}
\end{figure}

\begin{theorem} Graph Laplacian $L$ is nonnegative definite. 
\end{theorem}
The proof is based on factoring Laplacian $L$ using incidence matrix $\nabla$ such that
$L = \nabla^{\top} \nabla$. Such factorization always yields nonnegative definite matrices.
Very often $L$ is nonnegative definite in practice if it is too sparse (Figure \ref{fig:Laplacian-adj}).

\begin{theorem}For graph Laplacian $L$, $L + \alpha I$ is positive definite for any $\alpha>0$.
\end{theorem} 
{\em Proof.} Since $L$ is nonnegative definite, we have
$$ x^{\top} L x \geq 0.$$
Then it follows that
$$ x^{\top} (L + \alpha I) x =  x^{\top} L x + \alpha x^{\top} x > 0$$
for any $\alpha >0$ and $x \neq 0$. $\square$

\begin{figure}[t]
\centering
\includegraphics[width=0.8\linewidth]{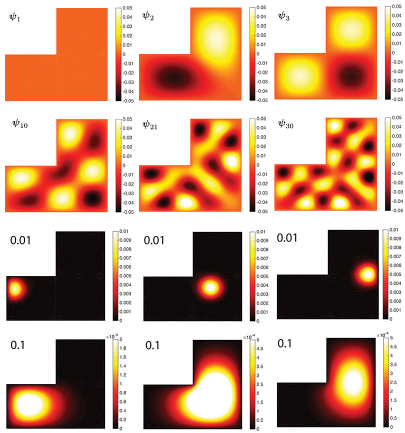}
\caption{Top: First few eigenvectors of the Laplacian in a L-shaped domain. Botom: Heat kernel with bandwidths $\sigma =0.01, 0.1$. We have used degree 70 expansions but the shape is almost identical if we use higher degree expansions. The heat kernel is a probability distribution that follows the shape of the L-shaped domain.}
\label{fig:Lshape-basis}
 \end{figure}

Unlike the continuous Laplace-Beltrami operators that may have possibly infinite number of eigenfunctions, we have up to $p$ number of eigenvectors $\psi_1, \psi_2, \cdots, \psi_p$ satisfying 
\bqn L \psi_j = \lambda_j \psi_j \label{eq:graph-eigen} \eqn
with (Figure \ref{fig:Lshape-basis})
$$0 = \lambda_1 < \lambda_2 \leq \cdots \leq \lambda_p.$$
The eigenvectors are orthonormal, i.e.,
$$\psi_i^{\top} \psi_j  = \delta_{ij},$$
the Kroneker's delta. The first eigenvector is trivially given as
$\psi_1 = {\bf 1}/\sqrt{p}$ with ${\bf 1} = (1, 1, \cdots, 1)^{\top}.$ 


All other higher order eigenvalues and eigenvectors are unknown analytically and have to be computed numerically (Figure \ref{fig:Lshape-basis}). Using the eigenvalues and eigenvectors, the graph Laplacian can be decomposed spectrally. From (\ref{eq:graph-eigen}), 
\bqn L \Psi = \Psi \Lambda, \label{eq:graph-psi}\eqn
where  $\Psi = [\psi_1, \cdots, \psi_p]$ and $\Lambda$ is the diagonal matrix with entries  $\lambda_1, \cdots, \lambda_p$. 
Since $\Psi$ is an orthogonal matrix, 
$$\Psi \Psi^{\top} = \Psi^{\top}\Psi = \sum_{j=1}^p \psi_j \psi_j^{\top} = I_p,$$ the identify matrix of size $p$. Then (\ref{eq:graph-psi}) is written as
$$L  = \Psi \Lambda \Psi^{\top} = \sum_{j=1}^p \lambda_j \psi_j \psi_j^{\top}.$$
This is the restatement of the singular value decomposition (SVD) for Laplacian. 

For measurement  vector $f = (f_1, \cdots, f_p)^{\top}$ observed at the $p$ nodes, the discrete Fourier series expansion is given by
$$f = \sum_{j=1}^n \tilde f_j \psi_j,$$
where $\tilde f_j =  f^{\top}\psi_j = \psi_j^{\top}f$ are Fourier coefficients.

 \section{Fiedler vectors} 
\index{Fiedler vector}

The connection between the eigenfunctions of continuous and discrete Laplacians have been well established by many authors \citep{gladwell.2002,tlusty.2007}. Many properties of eigenfunctions of  Laplace-Beltrami operator have discrete analogues. The second eigenfunction of the  graph Laplacian is called the Fiedler vector and it has been studied in connection to the graph and mesh manipulation, manifold learning and the minimum linear arrangement problem \citep{fiedler.1973,ham.2005,levy.2006,ham.2004,ham.2005}.

Let $G = \{V,E\}$ be the graph with the vertex set $V$ and the edge set $E$. We will simply index the node set as $V = \{1,2 ,\cdots, n\}$. If two nodes $i$ and $j$ form an edge, we denote it as $i \sim j$.  The edge weight between $i$ and $j$ is denoted as $w_{ij}$. For a measurement vector ${\bf f} = (f_1, \cdots, f_n)^{\top}$ observed at the $n$ nodes,  the discrete Dirichlet energy  is given by
\bqn \mathcal{E} ({\bf f}) = {\bf f}^{\top} L {\bf f} =  \sum_{i,j=1}^n w_{ij} (f_i - f_j)^2 = \sum_{i \sim j} w_{ij}(f_i - f_j)^2.\label{eq:DEdiscrete}\eqn
The discrete Dirichlet energy (\ref{eq:DEdiscrete}) is also called the linear placement cost in the minimum linear arrangement problem \citep{koren.2002}. Fielder vector ${\bf f}$ evaluated at $n$ nodes is obtained as the minimizer of the quadratic polynomial:
$$ \min_{\boldsymbol f}   \mathcal{E} (f)$$
subject to the quadratic constraint 
\bqn \| {\bf f} \|^2 = {\bf f}^{\top}{\bf f} = \sum_i f_i^2 = 1. \label{eq:fiedlerconstraint} \eqn 
The solution can be interpreted as the kernel principal components of a Gram matrix given by the generalized inverse of $L$ \citep{ham.2004,ham.2005}. 
Since the eigenvector $\psi_1$ of Laplacian is  orthonormal with eigenvector $\psi_0$, which is constant, 
we also have an additional constraint: 
\bqn \sum_{i} f_i =0.   \label{eq:fiedlerconstraint2} \eqn
This optimization problem was first introduced for the minimum linear arrangement problem in 1970's \citep{hall.1970,koren.2002}. The optimization can be solved using the Lagrange multiplier as follows \citep{holzrichter.1999}.
\index{graph!Laplacian}
\index{eigenfunctions!graph Laplacian}
\index{kernel!principal components}
\index{matrix!Gram}

Let $g$ be the constraint (\ref{eq:fiedlerconstraint}) so that
$$g({\bf f})  = {\bf f}^{\top}{\bf f} - 1 = 0.$$ Then the constrainted minimum should satisfy
\bqn \nabla \mathcal{E} - \mu \nabla g  =0, \label{eq:fiedler1}\eqn
where $\mu$ is the Lagrange multiplier. (\ref{eq:fiedler1}) can be written as
\bqn 2L{\bf f} - \mu {\bf f} =0 \label{eq:fiedler2}\eqn
Hence, ${\bf f}$ must be the eigenvector of $L$ and $\mu/2$ is the corresponding eigenvalue. By multiplying ${\bf f}^{\top}$ on the both sides of (\ref{eq:fiedler2}), we have
$$2{\bf f}^{\top}L{\bf f} = \mu {\bf f}^{\top} {\bf f} = \mu.$$
Since we are minimizing ${\bf f}^{\top} L {\bf f}$, 
$\mu/2$ should be the second eigenvalue $\lambda_1$.

In most literature \citep{holzrichter.1999}, the condition $\sum_{i} f_i = 0$ is incorrectly stated as a necessary constraint for the Fiedler vector. However, the constraint $\sum_i f_i = 0$ is not really needed in minimizing the Dirichlet energy. This can be further seen from introducing a new constraint 
$$ h({\bf f}) = {\bf e}^{\top}{\bf f} = \sum_i f_i = 0,$$
where ${\bf e} = (1, \cdots, 1)^{\top}$.

\begin{figure}[t]
\centering
\includegraphics[width=1\linewidth]{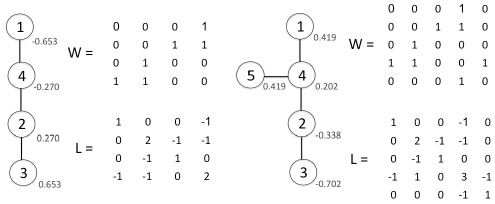}
\caption{A weighted graph with weights $W$ and the graph Laplacian $L$. The weights are simply the adjacency matrix. The second eigenvector $\psi_1$ is given as numbers beside nodes. Left: This example is given in \citet{hall.1970}. The maximum geodesic distance is obtained between the nodes 1 and 3, which are also hot and cold spots. Right: There are two hot spots 1 and 5 which corresponds to two maximal geodesic paths 1-4-2-3 and 5-4-2-3 \citep{chung.2011.MLMI}.}
\label{fig:graphlaplacian}
\end{figure}

The constraint (\ref{eq:fiedlerconstraint}) and (\ref{eq:fiedlerconstraint2}) forces $\psi_1$ to have at least two differing {\em sign domains} in which $\psi_1$ has one sign. But it is unclear how many differing sign domains $\psi_1$ can possibly have. The upper bound is given by Courant's nodal line theorem \citep{courant.1953,gladwell.2002,tlusty.2007}. The {\em nodal set} of eigenvector $\psi_i$ is defined as the zero level set $\psi_i(p) = 0$. Courant's nodal line theorem states that the nodal set of the $i$-th eigenvector $\psi_{i}$ divides the graph into no more than $i$ sign domains. Hence, the second eigenvector has exactly 2 disjoint sign domains. At the positive sign domain, we have the global maximum and at the negative sign domain, we have the global minimum. This property is illustrated in Figure \ref{fig:graphlaplacian}.
\index{Courant's nodal line theorem} 
However, it is unclear where the global maximum and minimum are located. The concept of tightness is useful in determining the location. 


\begin{definition}
For a function ${\bf f}$ defined on vertex set $V$ of $G$, let $G_s^-$ be the subgraph of $G$ induced by the vertex set $V_s^- = \{ i \in V  | f_i  < s \}$. Let $G_s^+$ be the subgraph of $G$ induced by the vertex set $V_s^+ = \{ i \in V | f_i  > s \}$. For any $s$, if $G_s^-$ and $G_s^+$ are either connected or empty, then ${\bf f}$ is tight  \citep{tlusty.2007}.
\end{definition}

\begin{figure}[t]
\centering
\includegraphics[width=0.8\linewidth]{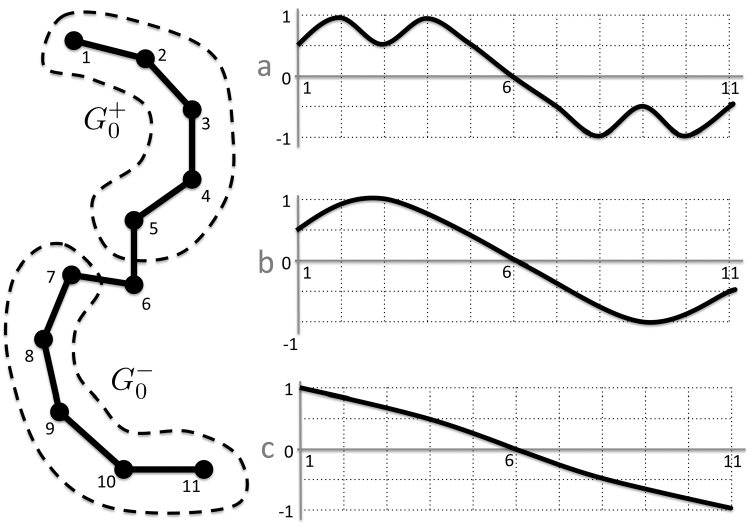}
\caption{A path with with positive ($G_0^+$) and negative ($G_0^-$) sign domains. Due to symmetry, the possible eigenfunction $\psi_1$ has to be an odd function.}
\label{fig:cycletight}
\end{figure}

When $s=0$, $G_0^+$ and $G_0^-$ are sign graphs. If we relax the condition so that $G_s^+$ contains nodes satisfying $f_i \geq s$, we have weak sign graphs. It can be shown that the second eigenvector on a graph with maximal degree 2 (cycle or path) is tight \citep{tlusty.2007}. Figure \ref{fig:cycletight} shows an example of a path with 11 nodes. Among three candidates for the second eigenfunction, (a) and (b) are not tight while (c) is. Note that the candidate function (a) have two disjoint components for $G_{0.5}^+$ so it can not be tight. In order to be tight,  the second eigenfunction cannot have a positive minimum or a negative maximum at the interior vertex in the graph \citep{gladwell.2002}. This implies that the second eigenfunction must decrease monotonically from the positive to negative sign domains as shown in (c). Therefore, the hot and cold spots must occur at the two end points $1$ and $11$, which gives the maximum geodesic distance of 11. 
\index{graph!degree}

For a cycle, the argument is similar except that a possible eigenfunction has to be periodic and tight, which forces the hot and cold spots to be located at the maximum distance apart. Due to the periodicity, we will have multiplicity of eigenvalues in the cycle.
Although it is difficult to predict the location of maximum and minimums in general, the behavior of the second eigenfunction is predictable for an elongated graph; it provides an intrinsic geometric way of establishing natural coordinates.

\section{Heat kernel smoothing on graphs}
\index{heat kernel}
\index{graph!heat kernel}

{\em Heat kernel smoothing} was originally introduced in the context of filtering out cortical surface data defined on mesh vertices obtained from 3D medical images  \citep{chung.2005.IPMI,chung.2005.NI}. The formulation uses the tangent space projection in approximating the heat kernel by iteratively applying Gaussian kernel with smaller bandwidth. Recently proposed spectral formulation to heat kernel smoothing \citep{chung.2015.MIA} constructs the heat kernel analytically using the eigenfunctions of the Laplace-Beltrami (LB) operator, avoiding the need for the linear approximation used in \citep{chung.2005.IPMI,han.2006}. Since surface meshes are graphs, heat kernel smoothing can be used to smooth noisy data defined on network nodes. 
		
Instead of Laplace-Beltrami operator for cortical surface, graph  Laplacian is used to construct the discrete version of heat kernel smoothing. The connection between the eigenfunctions of continuous and discrete Laplacians has been well established by several studies \citep{gladwell.2002,tlusty.2007}. Although many have introduced the discrete version of heat kernel in computer vision and machine learning, they mainly used the heat kernels to compute shape descriptors  or to define a multi-scale metric \citep{belkin.2006,sun.2009,bronstein.2010, deGoes.2008}. These studies did not use the heat kernel in filtering out  data on graphs. There have been significant developments in kernel methods in the machine learning community \citep{bernhard.2002,nilsson.2007,shawe.2004,steinke.2008,yger.2011}. However, the heat kernel has never been used in such frameworks. Most kernel methods in machine learning deal with the linear combination of kernels as a solution to penalized regressions. On the other hand, our kernel method does not have a penalized cost function.

\subsection{Heat kernel on graphs}
The {\em discrete heat kernel} $K_{\sigma}$ is a positive definite symmetric matrix of size $p \times p$ given by
\bqn K_{\sigma} = \sum_{j=1}^p e^{-\lambda_j \sigma} \psi_j\psi_j^{\top}, \label{eq:dhk} \eqn
where $\sigma$ is called the bandwidth of the kernel. Figure \ref{fig:Lshape-basis} displays heat kernel with different bandwidths at a L-shaped domain. Alternately, we can write (\ref{eq:dhk}) as
$$K_{\sigma} = \Psi e^{-{\sigma}\Lambda } \Psi^{\top},$$
where  $e^{-{\sigma}\Lambda }$ is the matrix logarithm of $\Lambda$. To see positive definiteness of the kernel,  for any nonzero $x \in \mathbb{R}^p$,
\bq
x^{\top} K_{\sigma} x &=& \sum_{j=1}^p e^{-\lambda_j \sigma} x^{\top}\psi_j\psi_j^{\top}x  \\
                &=& \sum_{j=1}^p e^{-\lambda_j \sigma} ( \psi_j^{\top}x )^2 > 0.\eq
When $\sigma=0$, $K_{0} = I_p$, identity matrix. When $\sigma=\infty$, by interchanging the sum and the limit, we obtain
$$K_{\infty} =  \psi_1\psi_1^{\top} = {\bf 1} {\bf 1}^{\top}/p.$$
$K_{\infty}$ is a degenerate case and the kernel is no longer positive definite. Other than these specific cases, the heat kernel is not analytically known in arbitrary graphs.

Heat kernel is doubly-stochastic \citep{chung.1997} so that 
$$K_{\sigma} {\bf 1} = {\bf 1}, \; {\bf 1}^{\top} K_{\sigma} = {\bf 1}^{\top}.$$
Thus, $K_{\sigma}$ is a probability distribution along columns or rows.

Just like the continuous counterpart,  the discrete heat kernel is also multiscale and has the scale-space property. Note
\bq K_{\sigma}^2 &=& \sum_{i,j=1}^p e^{-(\lambda_i + \lambda_j) \sigma} \psi_i\psi_i^{\top} \psi_j\psi_j^{\top}\\
&=& \sum_{j=1}^p e^{-2\lambda_j \sigma} \psi_j\psi_j^{\top} = K_{2\sigma}.
\eq
We used the orthonormality of eigenvectors. Subsequently, we have 
$$K_{\sigma}^n = K_{n\sigma}.$$

\subsection{Heat kernel smoothing on graphs}
\index{heat kernel smoothing}
\index{graph!smoothing}

Discrete heat kernel smoothing of measurement vector $f$  is then defined as convolution
\bqn K_{\sigma} * f = K_{\sigma}f  =\sum_{j=0}^p e^{-\lambda_j \sigma} {\tilde f}_j \psi_j, \label{eq:DHK} \eqn
This is the discrete analogue of heat kernel smoothing first defined in \citep{chung.2005.IPMI}. In discrete setting, the convolution $*$ is simply a matrix multiplication. Thus, $$K_{0} *f = f$$ and 
$$K_{\infty} *f =  \bar f {\bf 1},$$
where $\bar f = \sum_{j=1}^p f_j/p$ is the mean of signal $f$ over every nodes. When the bandwidth is zero, we are not smoothing data. As the bandwidth increases, the smoothed signal converges to the sample mean over all nodes.

Define the $l$-norm of a vector $f=(f_1, \cdots, f_p)^{\top}$ as 
$$\parallel f \parallel_{l} =
            \Big(  \sum_{j=1}^p  \big| f_j \big|^{l}  \Big)^{1/l}. 
$$
The matrix $\infty$-norm is defined as
$$
\parallel f \parallel_{\infty} =
            \max_{1 \leq j \leq p}  \big| f_j \big|. 
$$

\begin{theorem} Heat kernel smoothing is a contraction mapping with respect to the $l$-th norm, i.e.,
\bq \| K_{\sigma} *f \|_l^l \leq \|  f \|_l^l. \label{eq:contraction} \eq
\end{theorem}
{\em Proof.} Let kernel matrix $K_{\sigma} = (k_{ij})$. Then we have inequality
$$\| K_{\sigma} *f \|_l^l  = \sum_{i=1}^p  \sum_{j=1}^p |k_{ij}f_j |^l 
\leq \sum_{j=1}^p |f_j|^l.
$$
We used Jensen's inequality and doubly-stochastic property of heat kernel. Similarly, we can show that heat kernel smoothing is a contraction mapping with respect to the $\infty$-norm as well. 

Theorem 1 shows that heat kernel smoothing contracts the overall size of data. This fact can be used to skeltonize the blood vessel trees.

\subsection{Statistical properties}
Often observed noisy data $f$ on graphs is smoothed with heat kernel
$K_{\sigma}$ to increase the signal-to-noise ratio (SNR) and increases the statistical sensitivity
\citep{chung.2015.MIA}. We are interested in knowing how heat kernel smoothing will  affect on the statistical properties of smoothed data.

Consider the following addictive noise model:
\bqn f =    \mu + e, \label{eq:field} \eqn
where $\mu$ is unknown signal and $\epsilon$ is  zero mean noise.  Let $e= (e_1, \cdots, e_p)^{\top}$. Denote $\mathbb{E}$ as expectation and $\mathbb{V}$ as covariance. It is natural to assume that the noise variabilities  at different nodes are identical, i.e., 
\bqn \mathbb{E} e_1^2 = \mathbb{E} e_2^2 = \cdots = \mathbb{E} e_p^2. \label{eq:eeq} \eqn
Further, we assume that data at two nodes $i$ and $j$ to have less correlation when the distance between the nodes is large. 
So covariance matrix  $$R_{e}  =  \mathbb{V} e = \mathbb{E} (e e^{\top}) = (r_{ij})$$ can be given by
\bqn r_{ij} = \rho( d_{ij}) \label{eq:rij} \eqn
for some decreasing function $\rho$ and geodesic distance $d_{ij}$ between nodes $i$ and $j$. 
 Note
$r_{jj} = \rho(0)$ with the understanding that $d_{jj}=0$ for all $j$. The off-diagonal entries of $R_e$ are smaller than the diagonals.

Noise $e$ can be further modeled as Gaussian white noise, i.e., Brownian motion or  the generalized derivatives of Wiener
process, whose covariance matrix elements are Dirac-delta. For the discrete counterpart, $r_{ij} = \delta_{ij}$, where $\delta_{ij}$ is Kroneker-delta with
$\delta_{ij}=1$ if $i=j$ and 0 otherwise.  Thus, 
$$R_e = \mathbb{E} (e e^{\top}) = I_p,$$
the identity matrix of size $p \times p$. Since $\delta_{jj} \geq \delta_{ij}$, Gaussian white noise is a special case of (\ref{eq:rij}).

Once heat kernel smoothing is applied to (\ref{eq:field}), we have
\bqn K_{\sigma}*f = K_{\sigma}*\mu + K_{\sigma}*e. \label{eq:Kfield}\eqn
We are interested in knowing how the statistical properties of model change from (\ref{eq:field}) to (\ref{eq:Kfield}). 
For $R_e = I_p$, the covariance matrix of smoothed noise is simply given as
$$R_{K_{\sigma}*e}  = K_{\sigma} \mathbb{E} (e e^{\top}) K_{\sigma} =  K_{\sigma}^2 =K_{2\sigma}.$$
We used the scale-space property of heat kernel. In general, the covariance matrix  of smoothed data $K_{\sigma}* e$ is
given by
\bq R_{K_{\sigma}* e} = K_{\sigma} \mathbb{E} (e e^{\top}) K_{\sigma}  = K_{\sigma} R_e K_{\sigma}.\eq

The variance of data will be often reduced after heat kernel smoothing in the following sense \citep{chung.2005.IPMI,chung.2005.NI}: \\

\begin{figure}[t!]
\centering
\includegraphics[width=1\linewidth]{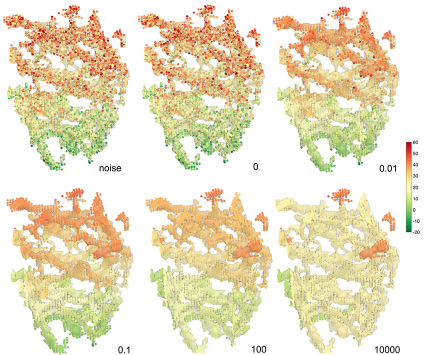}
\caption{From top left to right: 3D lung vessel tree. Gaussian noise is added to one of the coordinates. 3D graph constructed using 6-connected neighbors.  The numbers are the kernel bandwidth $\sigma$.} 
\label{fig:lung-tree}
\end{figure}

\begin{theorem}
\label{thm:diffusionvar}
Heat kernel smoothing reduces variability, i.e., 
$$\mathbb{V} (K_{\sigma}* f)_j \leq \mathbb{V} f_j \label{eq:varreduction}$$
for all $j$. The subscript $_j$ indicates the $j$-th element of the vector.
\end{theorem}
{\em Proof.} Note 
$$\mathbb{V} (K_{\sigma}* f)_j = \mathbb{V} (K_{\sigma}* e)_j = \mathbb{E} \Big( \sum_{i=1}^p k_{ij} e_i \Big)^2.$$
Since $(k_{ij})$ is doubly-stochastic, after applying Jensen's inequality, we obtain 
$$ \mathbb{E} \Big( \sum_{i=1}^p k_{ij} e_i \Big)^2 \leq \mathbb{E} \Big(  \sum_{i=1}^p k_{ij} e_i^2 \Big) = \mathbb{E} e_i^2.$$
For the last equality, we used the equality of noise variability (\ref{eq:eeq}). Since $\mathbb{E} f_j = \mathbb{E} e_i^2$, we proved the statement. $\square$

Theorem \ref{thm:diffusionvar} shows that the variability of data decreases after heat kernel smoothing. 

\subsection{Skeleton representation using heat kernel smoothing}

\begin{figure}[t]
\centering
\includegraphics[width=0.8\linewidth]{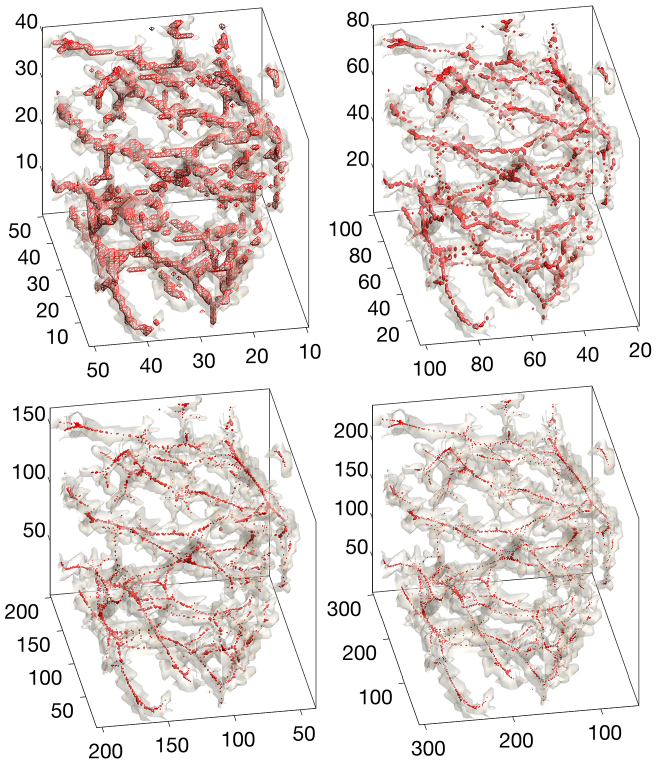}
\caption{The skeleton representation of vessel trees. Using the heat kernel series expansion with bandwidth $\sigma=1$ and 6000 basis, we upsampled the binary segmentation at 2, 4, 6 times (clockwise from top right) larger than the original size (top left).} 
\label{fig:lung-skel}
\end{figure}

Discrete heat kernel smoothing can be used to smooth out and present very complex patterns and get the skeleton representation. Here, we show how it is applied to the 3D graph obtained from the computed tomography (CT) of human lung vessel trees \citep{chung.2018.EMBCb,castillo.2009, wu.2013}. In this example, the 3D binary vessel segmentation from CT was obtained using the multiscale Hessian filters at each voxel \citep{frangi.1998,korfiatis.2011,shang.2011}. The binary segmentation was converted into a 3D graph by taking each voxel as a node and connecting neighboring voxels.  Using the 18-connected neighbor scheme, we connect two voxels only if they touch each other on their faces or edges.  If voxels are only touching at their corner vertices, they are not considered as connected. If the 6-connected neighbor scheme is used, we will obtain far sparse adjacency matrix and corresponding graph Laplacian. The eigenvector of graph Laplacian is obtained using  an Implicitly Restarted Arnoldi Iteration method \citep{lehoucq.1996}. We used 6000 eigenvectors. Note we cannot have more eigenvectors than the number of nodes.

As an illustration, we performed heat kernel smoothing on simulated data. Gaussian noise is added to one of the coordinates (Figure \ref{fig:lung-tree}). Heat kernel smoothing is performed on the noise added coordinate. Numbers in  Figure \ref{fig:lung-tree} are kernel bandwidths. At $\sigma=0$, heat kernel smoothing is equivalent to Fourier series expansion. Thus, we get the almost identical result. As the bandwidth increases, smoothing converges to the mean value.  Each disconnected regions should converge to their own different mean values.  Thus, when $\sigma=10000$,  the regions that are different colors are regions that are disconnected. This phenomena is related to the hot spots conjecture in differential geometry \citep{banuelos.1999, chung.2011.MLMI}.The number of disconnected structures can be obtained counting the zero eigenvalues.

The technique can be used to extract the skeleton representation of vessel trees. We perform heat kernel smoothing on node coordinates with $\sigma=1$. Then rounded off the smoothed coordinates to the nearest integers. The rounded off coordinates were used to reconstruct the binary segmentation. This gives the thick trees in Figure \ref{fig:lung-skel} (top left). To obtain thinner trees, the smoothed coordinates were scaled by the factor of 2, 4 and 6 times before rounding off. This had the effect of increasing the image size relative to the kernel bandwidth thus obtaining the skeleton representation of the complex blood vessel (Figure \ref{fig:lung-skel} clockwise from top right) \citep{lindvere.2013,cheng.2014}. By connecting the voxels sequentially, we can obtain the graph representation of the skeleton as well. The method can be easily adopted for obtaining the skeleton representation of complex brain network patterns.

\subsection{Diffusion wavelets}

Consider a traditional wavelet basis $W_{t,q}(p)$ obtained from a mother wavelet $W$ with scale and translation parameters $t$ and $q$ in Euclidean space \citep{kim.2012.NIPS}:
\bqn W_{t, q}(p) = \frac{1}{t}W \big(\frac{p-q}{t} \big). \label{eq:Wtq} \eqn
The wavelet transform of a signal $f(p)$ is given by kernel
\bq
     \langle W_{t,q}, f \rangle = \int_{\mathcal{M}} W_{t,q}(p)f(p) \;d\mu(p).
\eq

Scaling a function on an arbitrary manifold including graph is trivial. But the difficulty arises when one tries to translate a mother wavelet. It is not straightforward to generalize the Euclidean formulation (\ref{eq:Wtq}) to an arbitrary manifold, 
due to the lack of regular grids  \citep{nain.2007,bernal.2008}. The recent work based on the diffusion wavelet bypasses this problem also by taking bivariate kernel as a mother wavelet \citep{antoine.2010,hammond.2011,mahadevan.2006, kim.2012.NIPS}.  By simply changing the second argument of the kernel, it has the effect of translating the kernel. The diffusion wavelet construction has been fairly involving so far.  However, it can be shown to be a special case of the heat kernel regression with proper normalization. Following the notations in \citet{antoine.2010,hammond.2011,kim.2012.NIPS},  
 diffusion wavelet $W_{t,q}(p)$ at position $p$ and scale $t$ is given by
$$W_{t,q}(p) = \sum_{j=0}^k g(\lambda_j t)\psi_j(p)\psi_j(q),$$
for some scale function $g$. If we let $\tau_j = g(\lambda_j t)$, the diffusion wavelet transform  
is given by 
\bqn \langle W_{t,q}, f \rangle &=& \int_{\mathcal{M}} W_{t,q}(p)f(p) \;d\mu(p) \nonumber \\
& = &  \sum_{j=0}^k  g(\lambda_j t) \psi_j(q)   \int_{\mathcal{M}}  f(p) \psi_j(p) \; d\mu(p)  \nonumber \\
 &=& \sum_{j=0}^k \tau_j f_j \psi_j(q),   \label{eq:diffusion-waveletregression}
 \eqn
 where $f_j = \langle f, \psi_j \rangle$ is the Fourier coefficient. Note (\ref{eq:diffusion-waveletregression}) is the kernel regression \citep{chung.2015.MIA}. Hence, the diffusion wavelet transform can be simply obtained by doing the kernel regression without an additional wavelet machinery as done in \citet{kim.2012.NIPS}. 
Further, if we let $g(\lambda_j t) = e^{-\lambda_j t}$, we have 
$$W_{t,p}(q) = \sum_{j=0}^k e^{-\lambda_j t} \psi_j(p)\psi_j(q),$$
which is a heat kernel. The bandwidth $t$ of heat kernel controls resolution while the translation is done by shifting one argument in the kernel.

\section{Laplace equation} 
\index{cortical thickness!Laplace equation}
\index{Laplace equation}
\index{Glerkin's method}
\index{Iterative residual fitting}

\begin{figure}[t!]
\begin{center}
\includegraphics[width=1\linewidth]{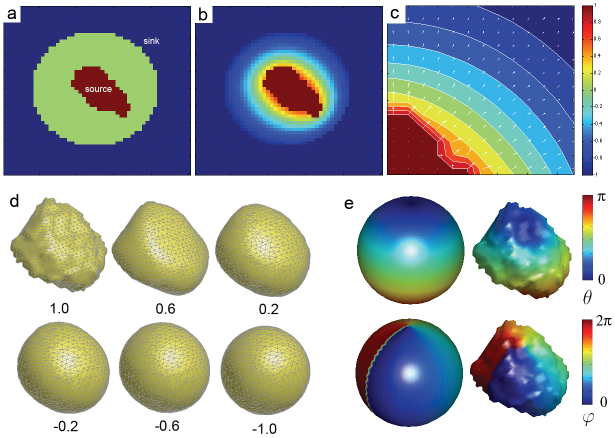}
 \caption{(a) The heat source (amygdala) is assigned value 1 while the heat sink (outer sphere) is assigned the value -1. The diffusion equation is solved with these boundary conditions. (b) After a sufficient number of iterations, the equilibrium state $f(\infty,p)$ is reached. (c) The gradient field $\nabla f(t=\infty,p)$ shows the direction of heat propagation from the source to the sink. The integral curve of the gradient field is computed by connecting one level set to the next level set of $f(\infty,p)$. (d) The deformation of amygala surface to the sphere is done by tracing the integral curve at each mesh vertex. The numbers $c=1.0, 0.6, \cdots, -1.0$ correspond to the level sets $f(\infty,p)$. (e) The surface flattening to the sphere produces the surface parameterization based on the spherical angles $(\theta, \varphi)$. The point $\theta=0$ corresponds to the north pole of a unit sphere. The method presented here is published in \citep{chung.2010.NI, chung.2012.CNA}.}\label{fig:surface-flattening}
\end{center}
\end{figure}

In this section, we will show how to solve for steady state of diffusion. The distribution of fictional charges within the two boundaries sets up a scalar potential field $f$, which satisfies the Poisson equation
$$\Delta f = \frac{\partial^2 f}{\partial x^2} + \frac{\partial^2 f}{\partial y^2} +\frac{\partial^2 f}{\partial z^2} = \frac{\rho}{\epsilon_0},$$
where $\rho$ is the total charge within the boundaries. If we set up the two boundaries at different potential, say at $f_{0}$ and $f_1$, without enclosing any charge, we have the Laplace equation
$$ \Delta f = 0.$$
The Laplace equation can be viewed as the steady state of diffusion
$$\frac{d f(t, p)}{d t} = \Delta f$$
when $ t \to \infty$. 
By solving the Laplace equation with the two boundary conditions, we obtain the potential field $f$. Then the electric field perpendicular to the isopotential surfaces is given by
$ - \nabla f$. The Laplace equation is mainly solved using the finite difference scheme \citep{chung.2012.CNA}. The electric field lines radiate from one conducting surface to the other without crossing each other. By tracing the electric field lines, we obtain one-to-one smooth map between surfaces (Figure \ref{fig:surface-flattening}). The underlying framework is identical to the Laplace equation based surface flattening or cortical thickness estimation \citep{jones.2000,chung.2010.NI}.

\subsection{Glerkin's method for solving Laplace equation}
 Without using the finite difference scheme, we can use an analytic approach for solving the Laplace equation using Galerkin's  method \citep{kirby.2000}. Galerkin's method discretizes partial differential equations and integral equations as a collection of linear equations involving  basis functions. The linear equations are then usually solved in the least squares fashion. The iterative residual fitting (IRF) algorithm \citep{chung.2007.tmi,chung.2008.sinica}, which iteratively fits functional data using diffusion equations, can be considered as a special case of Glerkin's method. 

We assume the eigenfunctions of $\Delta$ satisfying
$$\Delta \psi_j = \lambda_j \psi_j$$
are available. The solution of the Laplace equation is then approximated as a finite expansion
$$ f(p) = \sum_{j=0}^{k} c_j \psi_j(p).$$
Consider following boundary conditions
\bqn f(p) = 1, \; p \in G_+, \mbox{ and } f(p) = -1, \; p \in G_-, 
\eqn
where $G_+$ and $G_-$ are the subregions of $G$. $G_+$ and $G_1$ can be 3D objects, 2D surfaces, graphs or networks.  The integral curve between $G_1$ and $G_2$ will establish one to one correspondence.

The boundary conditions satisfy
\bqn 1 = \sum_{j=0}^k c_j \psi_j(p_{2i}), \;p_{2i} \in  G_+ \label{eq:galerkin1}\eqn
and
\bqn -1 = \sum_{j=0}^k c_j \psi_j(p_{3i}), \;p_{3i} \in G_- \label{eq:galerkin2}.\eqn
In the interior region $G \backslash ( G_+ \cup G_-)$, by taking the Laplacian on the expansion $f(p) = \sum_{j=0}^k c_j \psi_j (p)$, we have
\bqn 0 =  \sum_{j=0}^{k} c_j \Delta \psi_j(p_{1i})  = \sum_{j=0}^k c_j \lambda_j \psi_j(p_{1i}), \;p_{1i} \in G \backslash ( G_+ \cup G_-).\label{eq:galerkin3} \eqn

 We assume that there are $a, b$ and $c$ number of points for equations  (\ref{eq:galerkin1}), (\ref{eq:galerkin2}) and (\ref{eq:galerkin3}) respectively. We now combine linear equations (\ref{eq:galerkin1}), (\ref{eq:galerkin2}) and (\ref{eq:galerkin3}) together in a matrix form:
\bqn \underbrace{\left(\begin{array}{c}
0 \\
  \vdots \\
  0\\
 1\\
 \vdots\\
 1\\
 -1\\
 \vdots\\
 -1
\end{array}\right)}_{{\bf y}} =\underbrace{ \left (\begin{array}{ccc}
\lambda_0 \psi_0(p_{11})& \cdots & \lambda_k \psi_k(p_{11})\\
  \vdots & \ddots & \vdots \\
 \lambda_0 \psi_0(p_{1a})& \cdots & \lambda_k (p_{1a}) \\  
   \psi_0(p_{21}) & \cdots & \psi_k(p_{21}) \\
    \vdots & \ddots & \vdots \\
   \psi_0(p_{2b}) & \cdots & \psi_k(p_{2b}) \\
   \psi_0(p_{31}) & \cdots & \psi_k(p_{31}) \\
        \vdots & \ddots & \vdots \\
    \psi_0(p_{3c}) & \cdots & \psi_k(p_{3c}) \\
\end{array}\right)}_{\bf \Psi} \underbrace{\left(\begin{array}{c}
c_0 \\
c_1\\
  \vdots \\
  c_{k-1}\\
c_k
\end{array}\right)}_{\bf C}. \eqn
The size of matrix ${\bf \Psi}$ is $(a + b +c) \times (1+k)$. ${\bf \Psi} ' {\bf \Psi}$ is invertible if we sample substantially large number of points $ a + b +  c \gg k$. This is likely to be true in medical images so there is no need to use the pseudo-inverse here.  Then the matrix equation can be solved by the least squares method:
$$\widehat{\bf C} = ({\bf \Psi} ' {\bf \Psi})^{-1}{\bf \Psi} {\bf y}.$$

\subsection{Laplace equation with graphs}
We can solve for the Laplace equation within a graph $G$. We pick $G_+$ and $G_-$ to be subgraphs of $G$. 
For instance, we can take $G_+$ and $G_-$ at the two nodes in the graph and solve for the Laplace equation. With the boundary condition $f(G_+) = 1$ and $f(G_-)=1$, we are basically solving for steady state heat diffusion between heat source $G_+$ and heat sink $G_-$. 

Instead of solving the Laplace equation within a graph, we can also solve it between graphs by taking $G_+$ and $G_-$ to be two different graphs. For instance, we can take   two correlation matrices $G_+ =(g_{ij}^+)$ and $G_- = (g_{ij}^-)$, which can be viewed as weighted complete graphs. If we take $\Delta$ as the Hodge Laplacian defined on edges \citep{anand.2021.arXiv}, we can solve the Laplace equation with two correlation matrices as boundary conditions. This provides smooth one-to-one mapping between two correlation matrices. Unlike the usual element wise matching of $g_{ij}^+$ and $g_{ij}^-$, it provides the mapping through heat flow.

\section*{Acknowledgements}
The part of this study was supported by NIH grants NIH R01 EB022856 and R01 EB028753, NSF grant MDS-2010778. We would like to thank Yuan Wang of University of South Carolina for the providing Gibbs phenomenon plot and Gurong Wu  of University of Norh Carolina for providing lung CT scans. 

\bibliographystyle{agsm} 
\bibliography{reference.2021.12.28}

\end{document}